\documentclass[11pt]{amsart}
\usepackage{amssymb}
\numberwithin{equation}{section}

\newtheorem{theorema}{Theorem}
\newtheorem{theorem}{Theorem}[section]
\newtheorem{Theorem}[theorem]{Theorem}

\newtheorem{proposition}[theorem]{Proposition}
\newtheorem{lemma}[theorem]{Lemma}
\theoremstyle{definition}
\newtheorem{definition}[theorem]{Definition}
\theoremstyle{remark}
\newtheorem{remark}[theorem]{Remark}
\newtheorem{example}[theorem]{Example}

\newcommand\cS{\mathcal{S}}
\newcommand\cP{\mathcal{P}}
\newcommand\E{\mathcal{E}}
\newcommand{\U}{\mathcal{U}}

\newcommand{\Y}{\mathcal{Y}}
\renewcommand{\O}{\mathcal{O}}

\newcommand{\R}{\mathbb{R}}
\newcommand{\C}{\mathbb{C}}

\newcommand\lie[1]{\mathfrak{#1}}

\newcommand{\fh}{\lie{h}}
\newcommand{\fg}{\lie{g}}

\newcommand{\fk}{\lie{k}}

\def    \inv    {^{-1}}

\newcommand\cone{\mathaccent23c}

\newcommand{\SP}        {\operatorname{Sp}}

\begin{document}

\title{The Topological Structure of Contact and Symplectic Quotients }
\author{Eugene Lerman}
\author{Christopher Willett}
\address{Department of
Mathematics, University of Illinois, Urbana, IL 61801}
\email{lerman@math.uiuc.edu}
\email{cwillett@math.uiuc.edu}
\thanks{Partially supported by NSF grant DMS - 980305 and 
the American Institute of Mathematics.  C. Willett was 
supported by a graduate VIGRE fellowship.}

\begin{abstract}
We show that if a Lie group acts properly on a co-oriented contact
manifold preserving the contact structure then the contact quotient is
topologically a stratified space (in the sense that a neighborhood of
a point in the quotient is a product of a disk with a cone on a
stratified space). As a corollary we obtain that symplectic quotients
for proper Hamiltonian actions are topologically stratified spaces in
this strong sense thereby extending and simplifying the results of \cite{SL,
BL}.
\end{abstract}

\maketitle

\section{Introduction}
The notion of reduction of degrees of freedom in Hamiltonian mechanics
has a long history.  It was formalized in the early 1970s in the work of
Marsden and Weinstein and, independently, Meyer.  With the assumption
of a freely acting symmetry group, they showed that the reduced system
inherited a symplectic structure from the original system and that
the reduced dynamics was intimately related to the original dynamics.
A fair amount of effort by a number of people was spent during the
 1980s on
removing the assumption that the action is free.  For example, it was
proved by Arms, Marsden and Moncrief that the quotient of the zero
level set of a moment map by the group action is a union of symplectic
manifolds \cite{AMM}.  Arms, Cushman and Gotay showed that on
a symplectic reduced space there is 
a natural Poisson algebra of functions
\cite{ACG}.  In \cite{SL,BL} it was proved that reduced spaces are
stratified spaces in the sense that a neighborhood of a point is a
product of a disk with a cone on a compact 
stratified space, that the strata
are symplectic and that the symplectic structures on the strata are
tightly related.

The earliest appearance of the notion of a contact quotient that we
are aware of is in the work of Guillemin and Sternberg \cite{GS} on
homogeneous quantization, where it appeared in the guise of reduction
of symplectic cones.  Independently Albert \cite{?}
and Geiges \cite{Geiges} each showed
(under the assumption that the symmetry group acts freely on the zero
level set of the contact moment map) that the quotient of the zero
level set of the moment map by the 
group action, that is, the contact quotient, was
naturally a contact manifold.

The first result of this paper removes the freeness assumption and
shows that, in general, the contact quotients are naturally stratified
spaces:

 \begin{theorema}   \label{theoremA}
Let $M$ be a manifold with a co-oriented contact structure $\xi$.
Suppose a Lie group $G$ acts properly on $M$ preserving $\xi$.  Choose
a $G$-invariant contact form $\alpha$ with $\ker \alpha = \xi$ and let
$\Phi :M\to \fg^*$ be the corresponding moment map.

Then  for every subgroup $H$ of $G$, each connected component of
the topological space 
$$
\left( M_{(H)} \cap \Phi \inv (0) \right)/G
$$ 
is a manifold and the partition of the contact quotient 
$$
M_0 \equiv M/\!/G: = \Phi \inv (0)/G 
$$ 
into these manifolds is a stratification.  The symbol
$M_{(H)}$ stands for the set of points in $M$ with the isotropy groups
conjugate to $H$.
\end{theorema} 

It can be shown that the strata of a contact quotient are contact
manifolds.  This will be discussed elsewhere. 

The second main result of the paper is a natural short proof that
symplectic quotients are stratified spaces. The proof uses the 
fact that
neighborhoods of singularities of symplectic quotients are products of
a disk with a cone on a contact reduced space.  The result extends and simplifies the results of \cite{SL,BL}.

\begin{theorema}  \label{theoremB}
Let $(M, \omega)$ be a symplectic manifold with a proper Hamiltonian
action of a Lie group $G$ and a corresponding equivariant moment map
$\Phi: M\to \fg^* $.  Fix a point $\beta \in \fg^*$; denote its
isotropy group under the coadjoint action by $G_\beta$.

Then for every subgroup $H$ of $G_\beta$, each connected
component of the topological space 
$$
\left( M_{(H)} \cap \Phi \inv (\beta) \right)/G_\beta
$$ 
is a manifold, and the partition of the symplectic quotient at $\beta$
$$ 
M/\!/G (\beta) \equiv M_\beta :=
\Phi \inv (\beta)/G_\beta
$$ 
 into these manifolds is a stratification.  The symbol $M_{(H)}$ stands
for the set of points in $M$ with the isotropy groups conjugate to
$H$.
\end{theorema}

An outline of the contents of the paper is as follows:
\begin{enumerate}

\item Due to the multiple notions of stratification, we first make 
precise the meaning of the structure that we will place upon the 
contact and symplectic quotients.

\item A short review of the germane results of group actions on 
co-oriented contact manifolds is then pursued.  We define 
moment maps, the contact quotient, and give a sketch of the result of 
\cite{Geiges}.

\item We then prove the two main theorems simultaneously.    
The central idea is that both contact 
and symplectic quotients are locally modeled by reductions of vector
spaces by compact groups.

\item  We finish the paper with  a proof of a  local normal form theorem 
for group actions on contact manifolds which is modeled on the local
normal form theorem of Marle and of Guillemin and Sternberg for
symplectic manifolds.
\end{enumerate}

\subsection*{Acknowledgments}
The authors would like to thank Reyer Sjamaar for his helpful
suggestion that an inductive argument on dimension of the links of 
singularities should quickly  prove the main results of the paper.

\subsection*{A note on notation}  

Throughout the paper the Lie algebra of a Lie group denoted by a
capital Roman letter will be denoted by the same small letter in the
fraktur font: thus $\fg$ denotes the Lie algebra of a Lie group $G$
etc.  The identity element of a Lie group is denoted by 1.  The natural 
pairing between $\fg$ and $\fg^*$ will be denoted by 
$\langle \cdot, \cdot \rangle$.

When a Lie group $G$ acts on a manifold $M$ we denote the action by an
element $g\in G$ on a point $x\in G$ by $g\cdot x$ or occasionally by
$\tau_g(x)$; $G\cdot x$ denotes the $G$-orbit of $x$ and so on.  The
vector field induced on $M$ by an element $X$ of the Lie algebra $\fg$
of $G$ is denoted by $X_M$.  The isotropy group of a point $x\in M$ is
denoted by $G_x$; the Lie algebra of $G_x$ is denoted by $\fg_x$ and
is referred to as the isotropy Lie algebra of $x$.  We recall that
$\fg_x = \{ X \in \fg\mid X_M (x) = 0\}$.

If $P$ is a principal $G$-bundle then $[p, m]$ denotes the point in the
associated bundle $P\times _G M = (P\times M)/G$ which is the orbit of
$(p,m) \in P\times M$.

If $\omega$ is a differential form on a manifold $M$ and $Y$ is a 
vector field on $M$, the contraction of $\omega$ by $Y$ is 
denoted by $\iota(Y) \omega$.

\section{Preliminaries}
\subsection{Partitions and stratifications}

Since the word ``stratification'' is used to describe many different
types of partitions of topological spaces into manifolds, we start by
describing precisely the meaning that we attach to the term.  This is
a quick review.

\begin{definition}
A {\bf partition} of a Hausdorff topological space $X$ is a
collection of connected 
subsets $\{\cS_i\}$ of $X$ such that, set-theoretically,
$X$ is the disjoint union of $\cS_i$'s.  We  refer to the pair $(X,
\{\cS_i\})$ as a {\bf partitioned space} and to the subsets $\cS_i$ as {\bf
pieces}.
\end{definition}

\begin{remark}  The connectedness of the pieces in the definition 
of a partitioned space is slightly superfluous as we could simply 
refine a partition by non-connected pieces into one by connected 
pieces.
\end{remark}

\begin{remark}  
The product of two partitioned spaces $(X,\{\cS_i\})$ and 
$(Y, \{ \cP_i\})$ is the partitioned space $(X\times Y, \{\cS_i \times
\cP_j\})$.
\end{remark}

\begin{remark}  
An open subset $U$ of a partitioned space $(X,\{\cS_i\})$ is naturally
a partitioned space: the sets $U\cap \cS_i$ form a partition of $U$.
\end{remark}

\begin{definition}
Two partitioned spaces are {\bf isomorphic} if they are homeomorphic
and the homeomorphism takes pieces to pieces.
\end{definition}

Recall that a cone on a topological space $X$, denoted by 
$\cone(X)$,  is the quotient of the product $X\times [0,1)$ by the
relation $(x, 0) \sim (x', 0)$ for all $x, x'\in X$.  That is, $\cone
(X)$ is $X\times [0,1)$ with the ``boundary'' $X\times \{0\}$
collapsed to a point, the vertex $*$ of the cone.

\begin{remark}
The cone on a partitioned space $(X,\{\cS_i\})$ is the partitioned space
$$(\cone (X), \{*\} \bigsqcup\{\cS_i\times (0, 1)\})$$
\end{remark}

\begin{definition} 
If the pieces of a partitioned space $(X,\{\cS_i\})$ are
manifolds then we can define the dimension of $X$ to be 
$$
\dim X = \sup _i \dim \cS_i .
$$
\end{definition}

\begin{remark} Note that the dimension of a partitioned space 
depends upon the choice of the partition.  For example, every manifold
admits a partition into points. For such a partition the dimension of
the manifold as partitioned space is zero.  However, the manifold
dimension is always an upper bound for the partitioned space
dimension.
\end{remark}

We will only consider finite dimensional spaces.  Unless otherwise
mentioned we will partition connected manifolds into one piece only.
Note that if a space $X$ is partitioned into manifolds, then so is the
cone on $X$ and
$$
\dim \cone (X) = \dim X + 1.
$$

A stratification is a particularly nice type of a partition into
manifolds.  The definition is recursive on the dimension of
partitioned space.

\begin{definition} (cf.~\cite{GM}) 
A partitioned space $(X,\{\cS_i\})$ is a
{\bf stratified space} if each piece $\cS_i$ is a manifold and 
if for every piece $\cS$ and for every point
$x\in \cS$ there exist
\begin{enumerate}
\item an open neighborhood $U$ about $x$,
\item an open ball $B$ in $\cS$ about $x$,

\item a compact stratified space $L$, called the {\bf link} of 
the stratification of $X$ at $x$ ($L$ may be empty) and
 
\item an isomorphism  $\varphi : B \times \cone (L) \to U$
of partitioned spaces (If $L
 = \emptyset$ we require that $U$ is homeomorphic to the ball $B$.)
\end{enumerate}
The pieces of a stratified space are called {\bf strata}; the
collection of strata is called a {\bf stratification}.
\end{definition}

\begin{remark}
\begin{enumerate} 
\item A zero dimensional stratified space is a discrete set of points.

\item A connected manifold is a stratified space with exactly one stratum.

\item A cone on a compact manifold is a stratified space.
\end{enumerate}
\end{remark}

An example of a stratification is furnished by the orbit type
partition of an orbit space for a proper action of a Lie group on a
manifold.  Let us briefly review the relevant ideas (see the
introduction for the notational conventions).

An action of a Lie group $G$ on a manifold $M$ is {\bf proper} if the
map $G\times M \to M\times M$ given by $(g, m)\mapsto (g\cdot m, m)$
is proper, i.e., the preimages of compact sets under this map are
compact.  Note that for proper actions all isotropy groups are compact.

A {\bf slice} for an action of a Lie group $G$ on a manifold $M$ at a
point $x$ is a $G_x$ invariant submanifold $\Sigma$ such that $G\cdot
\Sigma$ is an open subset of $M$ and such that the map $G\times \Sigma
\to G\cdot \Sigma$, $(g, s) \mapsto g\cdot s$  descends to 
a diffeomorphism $G\times _{G_x} \Sigma \to G\cdot \Sigma$, $[g, s]
\mapsto g\cdot s$.  Thus for any point $y\in G\cdot \Sigma$, the 
orbit $G\cdot y$ intersects the slice $\Sigma$ in a single
$G_x$-orbit.  Also, for any $y\in \Sigma$, $G_y \subset G_x$.  A
theorem of Palais \cite{Palais} asserts that for smooth proper actions
slices exist at every point.

Let $M$ be a manifold with a proper action of a Lie group $G$.  For a
subgroup $H$ of $G$ denote by $M_{(H)}$ the set of points in $M$ whose
stabilizer is conjugate to $H$; i.e,
$$
M_{(H)}=\{x \in M\mid \text{there exists } g \in G \text{ such that } 
g G_x g\inv = H\},
$$
where $G_x$ denotes the stabilizer of $x$.  The set $M_{(H)}$ is often
referred to as the set of points of {\bf orbit type} $H$.  

It is an easy consequence of the existence of slices that the
connected components of the sets $M_{(H)}$ are manifolds and that the
components of the quotient $M_{(H)}/G$ are manifolds as well.
Moreover, using induction on the dimension of $M/G$ one can show that
the partition of $M/G$ into the connected components of the sets of
the form $M_{(H)}/G$ is a stratification.  The point of this paper is
to adapt the argument  for contact and symplectic quotients.

Note that in the paper we refer to both partitions $$M =
\bigsqcup_{(H)} M_{(H)}$$ and $$M/G = \bigsqcup_{(H)} M_{(H)}/G$$ as
stratifications by orbit type.

\begin{remark}  If $H$ and $K$ are conjugate subgroups of 
$G$, then $M_{(H)}=M_{(K)}$.  Thus, the indexing set for the 
stratification by orbit type is the set of conjugacy 
class of $G$.
\end{remark}

We will extensively use the following easy fact.

\begin{lemma}\label{easy_fact}
Suppose a Lie group $G$ acts properly on a manifold $M$.  Let $x\in M$
be a point and let $\Sigma \ni x$ be a slice through $x$ for the
action of $G$.  Denote the isotropy group of $x$ by $H$.

Then for any $H$-invariant subset $Z $ of $\Sigma$ 
$$
(G\cdot Z)/G = Z/H  \quad (\text{as partitioned spaces}) 
$$
where the left hand side is partitioned by $G$-orbit types and the
right hand  by $H$-orbit types.
\end{lemma}
\begin{proof}
Since $\Sigma$ is a slice,  for any $m\in Z \subset \Sigma$, $G\cdot m
\cap Z$ is a single $H$-orbit.
\end{proof}

\subsection{Group actions on contact manifolds}
Recall that a (co-oriented) {\bf contact structure} on a manifold $M$
of dimension $2n + 1$ is a distribution $\xi$ on $M$ given globally as
the kernel of a 1-form $\alpha$ such that $\alpha \wedge (d \alpha)^n
\neq 0$.  The form $\alpha$ is called a {\bf contact form}.
Whenever convenient we will refer to a pair $(M, \alpha)$ or to a pair
$(M, \xi)$ as a contact manifold.

An action of a Lie group $G$ on a manifold preserving a contact
structure $\xi$ is called a {\bf contact action}.

\begin{remark}
If a Lie group $G$ acts properly on a manifold and preserves a
co-oriented contact structure $\xi$, then there exists a $G$-invariant
1-form $\alpha$ with $\ker \alpha = \xi$.  If $G$ is compact, the
proof of this assertion is easy: since $\xi$ is co-oriented, there
exists by definition a 1-form $\alpha_0$ with $\ker \alpha_0 = \xi$.
If $\alpha_0$ is not $G$-invariant, average it over $G$.  If $G$ is not
compact, the argument is only slightly more complicated.  One adopts
Palais's proof of the existence of invariant Riemannian metrics on
manifolds with proper group actions \cite{Palais} to ``average''
$\alpha_0$.  See, for example, \cite{IJL} for details.
\end{remark}

\begin{definition}
Suppose a Lie group $G$ acts  on a manifold $M$ preserving a contact form
$\alpha$.  We define the {\bf contact moment map} $\Phi:M \to \fg^*$ by
$$
\langle \Phi(x), X \rangle=\alpha_x(X_M(x))
$$ 
for all $X$ in the Lie algebra $ \fg$ of $G$ and all $x\in M$.
The
contact moment map is $G-$equivariant, where as usual $G$ acts on
$\fg^*$ by the coadjoint action (see \cite{Geiges}).  Hence, the $G$ action
descends to an action on the zero level set and we define the {\bf
contact quotient} (or {\bf contact reduction}) of $M$ by $G$ to be 
the topological space
$$
M_0 \equiv M/\!/G:=\Phi\inv(0)/G .
$$ 
We will use the symbols $M/\!/G$ and $M_0$ interchangeably.

We define
the {\bf canonical partition} of the quotient $M/\!/G$ to be the
 the connected components of the sets of the form 
$$ 
        (M_{(H)} \cap \Phi \inv (0))/ G   
$$
for all conjugacy classes $(H)$ of $G$.
\end{definition}

Similarly symplectic quotients also have a canonical partition.

\begin{definition} \label{def_symp_quotient}
Let $(M, \omega)$ be a symplectic manifold with a proper Hamiltonian
action of a Lie group $G$ and a corresponding equivariant moment map
$\Phi: M\to \fg^* $.    We define the symplectic quotient of $M$ at a 
point $\beta \in \fg^*$ to be the topological space
$$
M/\!/G (\beta) \equiv M_\beta := \Phi\inv (\beta) /G_\beta ,
$$
where $G_\beta$ denotes the isotropy group of $\beta$ under the coadjoint
action.

Next we define the {\bf canonical partition} of the quotient $M_\beta$
to be the connected components of the sets of the form 
$$ 
(M_{(H)} \cap \Phi \inv (\beta))/ G_\beta 
$$ 
for all conjugacy classes $(H)$ of $G$.
\end{definition}

\begin{remark} 
Suppose a Lie group $G$ acts properly on a manifold $M$ preserving a
contact form $\alpha$. Albert and, independently, Geiges showed that
if $G$ acts freely on $\Phi\inv(0)$, then $0$ is a regular value of
the corresponding moment map $\Phi$, the contact quotient $M_0$ is a
manifold and $\alpha |_{\Phi\inv(0)}$ descends to a contact form
$\alpha_0$ on $M$ \cite{?,Geiges}. 
See also \cite{GS} where the
result was obtained earlier in a different form (cf.\ the introduction).
\end{remark}

\begin{remark} \label{rm2.16}
Recall that the {\bf symplectization} of a contact manifold $(M,
\alpha)$ is the symplectic manifold $(M \times \R, d(e^t \alpha))$.
If a Lie group $G$ acts on $M$ and preserves $\alpha$ then the trivial
extension of the action of $G$ to $M\times \R$ is Hamiltonian and a
corresponding moment map $\Psi:M\times \R \to \fg^*$ is related to the
contact moment map $\Phi: M\to \fg^*$ by the formula
$$
\Psi(m,t)=-e^t \Phi (m) \quad \text{for all} \quad (m,t) \in M\times \R.
$$ In particular, $\Psi\inv(0)=\Phi \inv(0) \times \R$.  One can use
this to give an alternative proof that, in the case of free proper actions,
contact quotients are contact manifolds.
\end{remark}

\begin{remark} 
Recall that the {\bf contactization} of an exact symplectic manifold
$(M,d \lambda)$ is the contact manifold $(M \times
\R,\lambda+dt)$.
If a Lie group $G$ acts on $M$ preserving the 1-form $\lambda$ then
the action is Hamiltonian with moment map $\Phi:M \to \fg^*$ given by
$\langle \Phi(m),X
\rangle=-\lambda_m(X_M(m))$ for all $X\in \fg$, all $m\in M$.  
The trivial extension of the action of $G$ to $M \times \R$
 preserves the contact form.  By definition the contact
moment map $\Psi:M \times \R \to \fg^*$ is given by $\langle
\Psi(m,t), X \rangle = -\lambda_m(X_M(m))$. Thus $\Psi$ factors
as $\Psi=\Phi \circ \text{pr}$, where $\text{pr}:M \times \R \to M$
is the projection.  It follows that, as spaces, 
$$
(M \times\R)/\!/G=M/\!/G (0)\times \R. 
$$
If $G$ acts freely on $\Phi \inv (0)$, it is not hard to show that the
symplectic quotient at zero $M/\!/G (0)$ is an exact symplectic
manifold and that the two sides are the same as contact manifolds.
\end{remark}

It will be particularly useful for us to consider the contactization of
symplectic vector spaces.  Suppose $(V, \omega_V)$ is a symplectic
vector space.  The Lie derivative of $\omega _V$ with respect to
the radial vector field $R(v) := v$ is $2\omega_V$, hence 
$$ 
d( \frac{1}{2}\iota (R)\omega_V) = \omega _V.  
$$ 

Therefore $\frac{1}{2}\iota (R)\omega_V + dt$ is a contact form on
$V\times \R$.  Note that if $\rho: H \to \SP (V, \omega_V)$ is a
symplectic representation of a Lie group $H$ then, since the action of
$H$ preserves the radial vector field $R$, the 1-form
$\frac{1}{2}\iota (R)\omega_V$ is $H$-invariant. Hence if we trivially
extended the action of $H$ to the contactization $V\times \R$, the
contact form $\frac{1}{2}\iota (R)\omega_V + dt$ is $H$ invariant as
well.  We will refer to the contact manifold $(V\times \R,
\frac{1}{2}\iota (R)\omega_V + dt)$ as a {\bf contact vector space} and as
{\bf the} contactization of $(V, \omega _V)$.

\section{Proof of the main result}

Our proof of the main result --- Theorems~\ref{theoremA} and
\ref{theoremB} --- has two ingredients.  The first one is an
observation that arbitrary contact quotients are modeled on the
contact quotients of contact vector spaces and that arbitrary
symplectic quotients are modeled on the symplectic quotients of
symplectic vector spaces.  The second one is Lemma~\ref{vs_lemma},
below, which describes the structure of contact and symplectic quotients
of vector spaces.

 The idea that a symplectic quotient for a reasonable group
action can be modeled on a symplectic quotient of a vector space is
due to R. Sjamaar.  The idea was implemented for reduction at zero by
an action of a compact Lie group in \cite{SL} and for reduction at a
locally closed coadjoint orbit by a proper action of an arbitrary Lie
group in \cite{BL}.

It was pointed out by Ortega and Ratiu that the technical assumption
on the coadjoint orbit in \cite{BL} is unnecessary if one defines
symplectic quotients as in Definition~\ref{def_symp_quotient}.  It was
shown in \cite{BL} that if the coadjoint orbit $G\cdot \beta$ is
locally closed then the sets $(M_{(H)} \cap \Phi \inv (G\cdot \beta))/
G \simeq (M_{(H)} \cap \Phi \inv (\beta))/ G_\beta $ are symplectic
manifolds and that the partition of $M_{G\cdot \beta} := \Phi \inv
(G\cdot \beta)/G \simeq M_\beta$ by these manifolds is a
stratification.  Ortega proved that one can drop the assumption on the
coadjoint orbit and still show that these sets are symplectic
manifolds, that they are locally closed in the symplectic quotient and
that the canonical partition is locally finite \cite{O} (see also
\cite{OR}).  Theorem~\ref{theoremA} asserts that the canonical
partition is in fact a stratification.  Our proof of it does not
directly use Ortega's result.  Rather, it uses the following two
propositions, which are variations on an argument in \cite{SL}.

\begin{proposition} \label{symplectic_model}
Let $(M, \omega)$ be a symplectic manifold with a proper Hamiltonian
action of a Lie group $G$ and a corresponding equivariant moment map
$\Phi: M\to \fg^* $ and let $M_\beta := \Phi \inv (\beta)/G$ be the
symplectic quotient at  $\beta \in \fg^*$.  For any point $x\in
M_\beta$ there exist a neighborhood $\U$ of $x$ in $M_\beta$, a
symplectic vector space $V$ with a symplectic action of a 
compact Lie group
$H$, a neighborhood $\U_0$ of the image of $0$ in the symplectic quotient
$V_0$ and a homeomorphism $\varphi :\U\to \U_0$ which maps the pieces of
the canonical partition of $M_\beta \cap \U$ into pieces of the
canonical partition of $V_0
\cap \U_0$.
\end{proposition}

\begin{proposition} \label{contact_model}
Suppose a Lie group $G$ acts properly on a manifold $M$ preserving a
contact form $\alpha$. Let $\Phi :M\to \fg^*$ be the
corresponding moment map.  For any point $x$ in the contact quotient
$M/\!/G := \Phi \inv (0)/G$ there exist a neighborhood $\U$ of $x$ in
$M/\!/G$, a contact vector space $V$ with a contact action of a
compact Lie
group $H$, a neighborhood $\U_0$ of the image of $0$ in a contact
quotient $V/\!/H$ and a homeomorphism $\varphi :\U\to \U_0$ which maps
the pieces of the canonical partition of $M_\beta \cap \U$ into pieces
of the canonical partition of $(V/\!/H) \cap \U_0$.
\end{proposition}
We postpone the proof of these two results for a later section and
proceed with the study of the quotients of vector spaces:

\begin{lemma} \label{vs_lemma}
Let $\rho : H \to \SP (V, \omega)$ be a symplectic representation of a
compact Lie group $H$.
\begin{enumerate}
\item  
The symplectic quotient at zero $V/\!/ H (0)$ is isomorphic, as a 
partitioned space, to the product of a symplectic vector space $U$ and
a cone on the quotient of a standard contact sphere $S$:
$$ 
V/\!/ H (0) = U \times \cone (S/\!/ H).  
$$
\item  

The contact quotient of the contactization of $V$ is isomorphic, as a
partitioned space, to the product of the contact vector space $U\times
\R$ and a cone on the quotient of the standard contact sphere $S$:
$$
(V\times \R)/\!/ H = (U\times \R) \times  \cone (S/\!/ H). 
$$
\end{enumerate}

\end{lemma}

\begin{proof}
Let $U$ be the vector subspace of $H$-fixed vectors in $V$; it is a
symplectic subspace of $V$.  Let $W$ denote the symplectic
perpendicular to $U$.  It is a symplectic $H$-invariant subspace as
well.  Let $R$ as before denote the radial vector field on $V$.  Since 
$d( \frac{1}{2}\iota (R)\omega) = \omega$, the moment map $\Phi_V$ for 
the action of $H$ on $V$ is given by the formula  
$$ 
\langle \Phi (v) , X \rangle = (\frac{1}{2}\iota (R)\omega) (X_V (v))
$$ 
for all $v\in V$ and all $X$ in the Lie algebra $\fh$ of $H$, where
$X_V$ denotes the linear vector field induced by $X$ on $V$.
Hence 
$\Phi \inv (0) = U \times \Phi _W \inv (0)$ where $\Phi _W $ is the
restriction of $\Phi$ to $W$.  Note that $\Phi _W$ is also a moment
map for the action of $H$ on $W$.

Since $H$ is compact there exists an $H$ invariant complex structure
$J$ on $W$ compatible with $\omega_W:= \omega |_W$.  Let $S$ denote
the unit sphere with respect to the inner product $\omega _W (\cdot, J
\cdot)$.  Since $\Phi _W$ is homogeneous,  
$\Phi _W\inv (0) = \cone (S \cap \Phi_W \inv (0))$.  Since the action
of $H$ is linear, $\Phi _W \inv (0)/H = U \times \cone (( S \cap
\Phi_W \inv (0))/H)$.  Note that the restriction of 
$\frac{1}{2}\iota (R)\omega_W $ to $S$ is a contact form and that the
restriction of $\Phi_W$ to $S$ is the corresponding contact moment
map.  Therefore the space $( S \cap\Phi_W \inv (0))/H)$ is the contact 
quotient $ S/\!/H$ and the first part of the lemma follows.
 
The second part of the lemma follows as well since the contact moment
map $\Psi$ on the contactization $V\times \R$ is related to the
symplectic moment map $\Phi$ on $V$ by $\Psi (v, t) = \Phi (v)$ for all
$(v,t) \in V\times \R$ (see Remark~\ref{rm2.16}).
\end{proof}

\begin{remark}
If in the Lemma above the representation $\rho$ is trivial, then
$V/\!/ H (0) = V$ and $(V\times \R)/\!/ H = V \times \R$.  In
particular, if $V = \{0\}$, then $(V\times \R)/\!/ H = \R$.
\end{remark}

\begin{proof}[Proof of Theorems~\ref{theoremA} and \ref{theoremB}]
Lemma~\ref{vs_lemma} together with Proposition~\ref{contact_model}
imply that a contact quotient is locally isomorphic as a partitioned
space to a product of an odd-dimensional disk and a cone on the
contact quotient of a standard contact sphere.  It follows that the
pieces of the canonical decomposition of a contact quotient are odd
dimensional manifolds and that they are locally closed.

Similarly Lemma~\ref{vs_lemma} together with
Proposition~\ref{symplectic_model} imply that a symplectic quotient is
locally isomorphic as a partitioned space to a product of an even
dimensional disk and a cone on the contact quotient of a standard
contact sphere.  It follows that the pieces of the canonical
decomposition of a symplectic quotient are even dimensional manifolds.

We next argue by induction on the dimension of contact quotients that
contact quotients are stratified spaces.  The smallest dimension that a
contact quotient can have is one.  If it is one dimensional then it is
a one dimensional manifold (which need not be connected).

Assume now that the dimension of our contact quotient $X$ is bigger
than one and that any contact quotient $Z$ with $\dim Z < \dim X$ is a
stratified space.  By Proposition~\ref{contact_model} a neighborhood
of a point in $X$ is isomorphic, as a partitioned space, to a
neighborhood of a point in the contact quotient of a contact vector
space.  By Lemma~\ref{vs_lemma}, the contact quotient of a contact
vector space is isomorphic to the product of a contact vector space
and a cone on the contact quotient of a sphere.  By induction, the
contact quotient of the sphere is a stratified space.  Therefore $X$
is a stratified space by definition.

Finally, by Proposition~\ref{symplectic_model} and Lemma~\ref{vs_lemma} a
symplectic quotient is isomorphic to a product of a disk and a cone on
the contact quotient of a contact sphere.  By the previous paragraph,
the contact quotient of the sphere is a stratified space.  Therefore a
symplectic quotient is also a stratified space.
\end{proof}

Observe that we have, in fact, proved the following.

\begin{Theorem}  Let $(X,\{\cS_\alpha\})$ be a partitioned 
space with the property that a neighborhood of every point 
is isomorphic (as partitioned spaces) to a product of 
a disk with the cone on the contact quotient of a 
standard contact sphere.  Then the partition $\{\cS_\alpha\}$ 
is  a stratification.
\end{Theorem}

We finish the section with a short proof that connected contact and
symplectic quotients have unique connected open dense strata.  In
particular we recover Theorem~5.9 of \cite{SL}.  The proof has two
parts.  We first remark that the contact quotient of a standard
contact sphere by a linear action of a compact Lie group is connected.

\begin{lemma}\label{con_link}
Let $S^{2n-1} \subset \C^n$ be the  standard contact sphere and let $K$
be a closed subgroup of the unitary group $U(n)$.  Then the contact
quotient $S^{2n-1}/\!/ K$ is connected.
\end{lemma}
\begin{proof}
The map $f: S^{2n-1} \to \R$, $f(z) = \frac{1}{2} ||z||^2$ is a moment
map for the action $(\lambda, z)\mapsto \lambda z$ of the circle $S^1
=U(1)$ on $\C^n$, and the sphere $S^{2n-1}$ is a level set of $f$.
Since the action of $S^1$ commutes with the action of $K$, the
restriction of the $K$-moment map $\Phi :\C^n \to \fk^*$ to
$S^{2n-1}$ descends to a moment map $\bar{\Phi}$ for the action of $K$
on the symplectic quotient $f\inv (2)/S^1 = S^{2n-1} /S^1 = \C
P^{n-1}$.  Since the projective space $ \C P^{n-1}$ is compact and
connected, the fibers of the moment map $\bar{\Phi}: \C P^{n-1} \to
\fk^*$ are connected by a theorem of Kirwan (see
\cite[Remark~2.1]{Kconv} and \cite[Remark~9.1]{Kbook}).

On the other hand, the restriction $\Phi |_{S^{2n-1}}$ is the contact
moment map for the action of $K$ on the sphere $S^{2n-1}$.
Since the fibers of $\bar{\Phi}$ are $S^1$ quotients of the fibers of
$\Phi |_{S^{2n-1}}$, the connectedness of the fibers of $\bar{\Phi}$
implies that the fibers of$\Phi |_{S^{2n-1}}$ are connected as well.
Therefore $S^{2n-1}/\!/ K = (\Phi |_{S^{2n-1}})\inv (0)/K$ is connected.
\end{proof}

The second part of the proof is:
\begin{proposition}  \label{open_dense} Let $(X,\{\cS_i\})$ be a connected 
stratified space. Assume, recursively, that all links are connected.
That is, the link of every point is connected, the link of every point
in the link is connected and so on.  Then there is a unique open dense 
stratum $X^r$ in $X$. 
\end{proposition}

\begin{proof}  The proof is an induction on the dimension of $X$. Let $X^r$ be 
the union of all the open strata in $X$.  We show first that $X^r$ is
dense.  Using density, we then show that $X^r$ is connected and hence
consists of a single stratum.  Note that a point in a stratified space
has an empty link if and only if it lies in an open stratum.  This
implies that if a stratum contains a set which is open in $X$ then the
whole stratum is open as well.

Let $x$ be a point in $X$ and $ \cS$ be the stratum containing $x$.
  By definition there an open neighborhood $U$ of $x$ in $X$,
an open ball $B$ in $\cS$, and a isomorphism $\rho:\cone(L)
\times B \to U$ of partitioned spaces, where $L$ is the link of $x$.
If the link is empty, then $x \in X^r$.  Otherwise, by induction $L$
contains a unique open dense stratum $L^r$. Then $L^r \times
(0,1)\times B$ is open in $\cone(L)\times B$.  Let $\cS'$ denote the
stratum in $X$ with $\rho (L^r \times (0,1)\times B) = \cS' \cap U$.
Then $ \cS' \cap U$ is open in $U$.  Therefore $\cS'$ is open in $X$,
hence $\cS' \subset X^r$.  Clearly $x$ lies in the closure of $\cS'$.
Therefore $x$ is in the closure of $X^r$.  This proves that $X^r$ is
dense.

We now prove that $X^r$ is connected.  Suppose not. Choose disjoint
open sets $U_1$ and $U_2$ in $X$ with $U_1 \cup U_2=X^r$.  By density,
$X=\bar U_1 \cup \bar U_2$.  Because $X$ is connected, however, $\bar
U_1 \cap \bar U_2 \neq
\emptyset$.
Choose $y \in \bar U_1 \cap \bar U_2$.  Then $y \not \in X^r$ and so
the link $L_y$ of $y$ is non-empty.  Hence, by induction, there is a
unique open stratum ${L_y}^r$ in $L_y$.  Choose a neighborhood, $U$ of
$y$ of the form $U= \cone(L_y)\times B$, where $B$ is a ball.  On the
one hand, $X^r \cap U= {L_y}^r \times (0,\epsilon ) \times B$ is
connected.  On the other, $X^r \cap U=(U_1 \cap U) \cap (U_2 \cap U)$
is not.
\end{proof}

\begin{remark}  We can remove the hypothesis that the space $X$ is connected. 
In this case, we work component by component to conclude that each
connected component of $X$ has a unique connected open dense stratum.
\end{remark}

Putting Lemma~\ref{con_link} and Proposition~\ref{open_dense} together
and using the proof of Theorems~1 and 2, we obtain:
\begin{theorem}  Suppose a Lie group $G$ acts properly on a
manifold $M$ preserving a contact form $\alpha$.  Then each connected
component of the contact quotient $M/\!/G$ has a unique connected open dense
stratum.

Let $(M, \omega)$ be a symplectic manifold with a proper Hamiltonian
action of a Lie group $G$ and a corresponding equivariant moment map
$\Phi: M\to \fg^* $.  Fix a point $\beta \in \fg^*$; denote its
isotropy group under the coadjoint action by $G_\beta$.  Then for
every point $\beta \in \fg^*$ a connected component of the symplectic
quotient $M/\!/G (\beta)$ has a unique connected open dense stratum.
\end{theorem}

\section{Local models of symplectic and contact quotients}

The main tool for proving Proposition~\ref{contact_model} is 

\begin{theorem} \label{contact_normal_form}
Suppose a Lie group $G$ acts properly on a manifold $M$ preserving a
contact form $\alpha$ and let
  $\Phi :M\to
\fg^*$ be the corresponding moment map.  Let $x\in \Phi \inv (0)$ be a point, let $H$ be its isotropy
group.  Denote the Lie algebra of $H$ by $\fh$ and the annihilator of
$\fh$ in $\fg^*$ by $\fh^\circ$.  Choose an $H$-equivariant splitting
$\fg^* =
\fh^\circ \oplus \fh^*$; let $i:\fh^* \to \fg^*$ be the corresponding
injection.

There exists a $G$-invariant neighborhood $U$ of $x$ in $M$, a
$G$-invariant neighborhood $U_0$ of the zero section on the vector
bundle $\Y = (G\times_H (\fh^\circ \times W))  \to G/H$ and a
$G$-equivariant diffeomorphism $\phi: U_0 \to U$ such that 
$\phi([1, 0, 0]) = x$ and
$$
(\Phi \circ \phi) ([g, \eta, w])
    = f([g, \eta, v])\, Ad^\dagger (g) \left(\eta + i (\Phi_W (w))\right),
$$
where
\begin{enumerate}
\item  $Ad^\dagger : G \to GL (\fg^*)$ is the coadjoint representation,

\item 
$W$ is the contactization of the maximal symplectic subspace $V$ of the
symplectic vector space $(\ker \alpha _x, d\alpha _x)$ complementary
to the tangent space $T_x (G\cdot x)$; $V$ can be and is chosen to be
$H$-invariant;

\item $\Phi_W : W \to \fh^*$ is the moment map for the linear action of $H$ on 
the contact vector space $W$;

\item $ f$ is a nowhere vanishing function;

\item $G$ acts on $\Y$ by $g\cdot [a, \eta ,v] = [ga, \eta, v]$.
\end{enumerate}
\end{theorem} 
We postpone the proof of this Theorem till the last section of the
paper and proceed with the proof of Proposition~\ref{contact_model}.

\begin{proof} [Proof of Proposition~\ref{contact_model}]

Define $F: \Y \to \fg^*$ by $F([g, \eta, w]) = Ad^\dagger (g)
\left(\eta + i (\Phi_W (w))\right)$. Since the function $f$ is
nonvanishing the map $\phi\inv$ sends $\Phi \inv (0) \cap U$ equivariantly
and homeomorphicly onto the set $F\inv (0) \cap U_0$.  Hence $(\Phi
\inv (0) \cap U)/G = (F\inv (0) \cap U_0)/G$ as partitioned spaces
(the partitioning is by $G$-orbit type).  Therefore it is enough to
show that
\begin{equation}\label{eq1}
F\inv (0)/G = \Phi_W \inv (0)/H \equiv W/\!/H \quad 
\text{as partitioned spaces.}
\end{equation}
Next note that the vector space $\fh^\circ \times W$ embeds
canonically into $\Y$ by $(\eta, w) \mapsto [1, \eta, w]$ and that
$\fh^\circ \times W$ is a slice through $[1, 0, 0]$ for the action of
$G$ on $\Y$.

Because $\fh^\circ \oplus i (\fh^*) = \fg^*$ and $i:\fh^* \to
\fg^*$ is injective, $\eta + i (\Phi _W (w)) = 0$ if and only if $\eta
= 0$ and $\Phi _W (w) = 0$.  Therefore 
$$ 
F\inv (0) = \{[g, \eta ,w]
\mid \eta = 0 \quad \text{and } \quad \Phi _W (w) = 0\} = G \times _H
(\{0\} \times \Phi _W \inv (0)).  
$$ 
Hence $F\inv (0) = G\cdot \Phi _W \inv (0)$.  Therefore by 
Lemma~\ref{easy_fact} equation~(\ref{eq1}) holds.
\end{proof}

\begin{proof}[Proof of Proposition~\ref{symplectic_model}]

The proof is essentially the same as that of
Proposition~\ref{contact_model} above.  Since we are not reducing as
zero, it is a little more delicate.  We use a version of the local
normal form theorem due to Marle and to Guillemin and Sternberg which
is proved on pp. 212 -- 215 of \cite{BL}:

\begin{theorem} \label{symplectic_normal_form}
Let $(M, \omega)$ be a symplectic manifold with a proper Hamiltonian
action of a Lie group $G$ and a corresponding equivariant moment map
$\Phi: M\to \fg^* $.  Let $x\in M$ be a point. Let $H$ denote its
isotropy group. Let $\beta = \Phi (x)$ and let $G_\beta$ denote the
isotropy group of $\beta$.  Define
$$ 
V := T_x (G\cdot x)^\omega / ( T_x (G\cdot x)^\omega \cap T_x (G\cdot x)), 
$$ 
where $T_x (G\cdot x)^\omega$ denotes the symplectic perpendicular to
the tangent space to the orbit $T_x (G\cdot x)$ in $T_x M$; it is a
symplectic vector space with a linear symplectic action of $H$.

Let $\fg$, $\fg_\beta$ and $\fh$ denote the Lie algebras of $G$,
$G_\beta$ and $H$, respectively. Let $\fh^\circ$ denote the
annihilator of $\fh$ in $ \fg_\beta^*$, and
$\fg_\beta ^\circ$ the annihilator of $\fg_\beta$ in $\fg^*$.  Choose
an $H$-equivariant splitting $\fg^* = \fh^* \oplus \fh^\circ \oplus
\fg_\beta ^\circ$; let $i : \fh^* \to \fg_\beta ^* $ and 
$j: \fg_\beta^* \to \fg^*$ denote the corresponding injections.

There exists a $G$-invariant neighborhood $U$ of $x$ in $M$, a
$G$-invariant neighborhood $U_0$ of the zero section on the vector bundle
$\Y = (G\times_H (\fh^\circ \times V)) \to G/H$ and a $G$-equivariant
diffeomorphism $\phi: U \to U_0$ such that
$$
(\Phi \circ \phi) ([g, \eta, v]) 
         =  Ad^\dagger (g) \left(\, \beta + j(\eta + i(\Phi_V (v))) \,\right),
$$
where
\begin{enumerate}
\item  $Ad^\dagger : G \to GL (\fg^*)$ is the coadjoint representation and 

\item $\Phi_V : V \to \fh^*$ is the moment map for the  action of $H$ on 
 $V$.

\end{enumerate} 
\end{theorem}
Define $F: \Y \to \fg^*$ by $F([g, \eta, w]) = Ad^\dagger (g)\left(\,
\beta + j(\eta + i(\Phi_V (v))) \,\right)$.  The map $\phi$ sends
$\Phi \inv (\beta)\cap U$ homeomorphicly and $G$-equivariantly onto
$F\inv (\beta) \cap U_0$.  Hence $(\Phi\inv (0) \cap U)/G_\beta =
(F\inv (0) \cap U_0)/G_\beta$ as spaces partitioned by $G$-orbit
types.  Therefore it is enough to prove that for some sufficiently
small $G_\beta$-invariant neighborhood $\O$ of $[1,0,0]$ in $\Y$ 
\begin{equation} \label{eq2.0}
(F\inv (0) \cap \O)/G_\beta = (\Phi _V \inv (0) \cap \tilde{\O})/H 
\end{equation}
as partitioned spaces for some $H$-invariant neighborhood $\tilde{\O}$
of 0 in $V$, where the left hand side is partitioned by $G$-orbit
types and the right hand side is partitioned by $H$-orbit types.

Note that the canonical embedding of the vector space $\fh^\circ
\times V$ into $\Y$, $(\eta , v) \mapsto [1, \eta, v]$ makes $\fh^\circ
\times V$ into a slice at $[1, 0, 0]$ for the action of $G$ on $\Y$.
The vector space $\fh^\circ \times V$ is also a slice at $[1, 0, 0]$
for the action of $G_\beta$ on $G_\beta \times _H ( \fh^\circ \times
V)$.  Therefore $\Y/G =(\fh^\circ
\times V)/H = \left(G_\beta \times _H ( \fh^\circ \times
V)\right) /G_\beta$.  It follows from Lemma~\ref{easy_fact} that for
any subgroup $K$ of $G_\beta$
\begin{equation} \label{eq2}
(G_\beta \times _H ( \fh^\circ \times V))_{(K)} = 
Y_{(K)} \cap \left(G_\beta \times _H ( \fh^\circ \times V)\right).
\end{equation}
Conversely if $Y_{(K)} \cap \left(G_\beta \times _H ( \fh^\circ \times
V)\right) \not = \emptyset$, we may assume that $K\subset
G_\beta$. Then equation (\ref{eq2}) holds.  Therefore the partition of
the left hand side of (\ref{eq2.0}) by $G$-orbit types is the same as
its partition by $G_\beta$-orbit types.

Recall that the tangent space to the coadjoint orbit through $\beta$
is canonically isomorphic to the annihilator of $\fg_\beta$ in
$\fg^*$: $T_\beta (Ad^\dagger (G) \beta) \simeq \fg_\beta ^\circ$.
Since $\fg^* = j (\fg_\beta^*) \oplus \fg_\beta ^\circ$, the vector
bundle $G\times _{G_\beta} \fg_\beta^*$ is the normal bundle of the
orbit $Ad^\dagger (G) \beta$ in $\fg^*$.  Consider the map $\E
:G\times _{G_\beta} \fg_\beta^* \to \fg^*$ given by $\E ([g, \eta])
=Ad^\dagger (g)(\beta + j (\eta))$.  It is the exponential map for a
flat $H$-invariant metric on $\fg^*$.  The differential $d\E$ is an
isomorphism at the point $[1, 0]$.  Therefore $\E$ is an open
embedding on a sufficiently small neighborhood $\O'$ of $[1, 0]$.  In
particular $\O'\cap \E \inv (\beta) = \{[1,0]\}$.  We now factor $F$ as
follows: let 
$$ 
F_2: G\times _H (\fh^\circ \times V) \to G\times _H
\fg_ \beta ^*, \quad 
[g, \eta , v] \mapsto [g, \eta + i ( \Phi _W (v))],
$$
 
$$
F_1: G\times _H \fg_ \beta ^* \to G\times _{G_\beta} \fg_ \beta ^*, 
\quad 
[g, \eta ] \mapsto [g, \eta ].
$$

We have $F = \E \circ F_1 \circ F_2$.  Let $\O = (F_1 \circ F_2)\inv
(\O ')$. Then $F\inv (\beta) \cap \O = (F_1 \circ F_2)\inv (\O ' \cap
\E \inv (\beta)) = \left( G_\beta \times _H (\{0\} \times \Phi _V \inv
(0)) \right) \cap \O$.  Consequently by Lemma~\ref{easy_fact} $(F\inv
(\beta) \cap \O) /G_\beta = (\Phi _V \inv (0) \cap \tilde{ \O})/H$ as
partitioned spaces, where the left hand side is partitioned by
$G_\beta$-orbit types, the right hand side is partitioned by $H$-orbit
types and where $ \tilde{ \O} = V \cap \O$.  But we have seen that the
partition of $G_\beta \times _H (\fh^\circ \times V)$ by $G$-orbit
types and by $G_\beta$ orbit types is the same partition.  The
Proposition now follows.
\end{proof} 

\section{Contact local normal form}\label{sec_normal_form}

In this section we provide a proof of the local normal form theorem 
(Theorem ~\ref{contact_normal_form}) that allowed us to establish 
our main results (Theorems ~\ref{theoremA} and \ref{theoremB}).
Not surprising, this requires an equivariant Darboux theorem for 
contact manifolds, the proof of which is the standard deformation 
argument of Moser.  We then review the contact analogues of 
isotropic submanifolds and symplectic normal bundles, introduced by 
Weinstein in \cite{wein}.  We finish with the proof of the local 
normal form theorem.

Recall that the Reeb vector field $Y$ on a contact manifold $(M,
\alpha)$ is the unique vector field satisfying $\iota(Y)\alpha=1$ and
$\iota(Y) d \alpha=0$.  If a Lie group acts on $M$ preserving
$\alpha$ then by uniqueness the Reeb vector field is invariant.

Observe that an equivariant version of Gray's stability theorem holds.

\begin{Theorem} (Equivariant Gray's Stability Theorem) 
Suppose a Lie group $G$ acts on a contact manifold $M$ and
$\{\alpha_t\}, t \in [0,1]$ is a smoothly varying
family of invariant
 contact forms. Suppose that $N \subseteq M$ is a invariant
submanifold and that
$$
\frac{\text{d}}{\text{dt}} \alpha_t(x)=
        \frac{\text{d}}{\text{dt}} d \alpha_t(x)=0
$$ 
for all $x \in N$. Then there is a family of equivariant
diffeomorphisms $\{\psi_t\}$ of $M$ with $\psi_{t |N}=\text{Id}$ and
$\psi_t^*(\alpha_t)=f_t \alpha_0$ for all $t$ and for some family of
positive functions $\{f_t\}$.
\end{Theorem}

As a consequence of Gray's stability theorem we get:

\begin{Theorem}\label{ercd} (Equivariant Relative Contact Darboux) 
Let a Lie group $G$ act properly on a manifold $M$ preserving an
embedded submanifold $N$.  Suppose $\alpha$ and $\beta$ are two
invariant contact forms with
$$\aligned
\alpha(x) &= \beta(x) \cr
d \alpha(x) &= d \beta(x) \endaligned
$$
for all $x \in N$.
Then there are invariant open neighborhoods $U_0,U_1$ of $N$, an equivariant 
diffeomorphism $\psi:U_0 \to U_1$ fixing $N$, 
and a positive function $f$ such that $\psi^*(\beta)=f \alpha$.
\end{Theorem}

\begin{proof} For $t \in [0,1]$, define $\gamma_t=(1-t) \alpha + t \beta$.
For all $x \in N$ and all $t$, we have
$\gamma_t(x)=\alpha(x)=\beta(x)$ and $d \gamma_t(x)=d \alpha(x)=d
\beta(x)$.  Hence, $\gamma_t$ is contact in an open neighborhood of
$N$.  Since $\gamma_t$ is invariant, Gray's theorem applies.  The time
1 map of the isotopy $\psi_t$ exists on some open set $O$ of $N$ since
it exists on a neighborhood of each point of $N$ in $M$.  Set $U_0=G
\cdot O$ and $U_1=\psi_1(U_0)$.
\end{proof}

\begin{definition} Let $(M,\alpha)$ be a 
contact manifold.  A submanifold $L$ of $M$ is {\em isotropic} if it
is tangent to the contact distribution $\ker \alpha$, that is, $T_xL \subset
\ker \alpha _x$ for every $x\in L $.
\end{definition}

\begin{example}  Suppose a Lie group $G$ acts properly on a 
 manifold $M$ preserving a contact form $\alpha$.  Then orbits $G
 \cdot x$ which lie in the zero level set of the contact moment
 $\Phi$ are isotropic submanifolds by definition of the moment map.
\end{example}

\begin{remark} \label{rm_isotr}
If $L$ is an isotropic submanifold of a contact manifold $(M,\alpha)$
then for every point $x\in L$ the tangent space $T_xL$ is an isotropic
subspace of the symplectic vector space $(\ker \alpha_x, d\alpha_x)$,
hence the terminology.
\end{remark}

\begin{definition}
Let $L$ be an isotropic submanifold of a contact manifold
$(M,\alpha)$.  Then for every point $x\in L$ the quotient
$(T_xL)^{d\alpha}/T_x L$ is a symplectic vector space, where
$(T_xL)^{d\alpha}$ denotes the symplectic perpendicular to $T_xL$ in
$\ker \alpha _x$ with respect to $d\alpha$ ( $T_x L \subset
(T_xL)^{d\alpha}$ by Remark~\ref{rm_isotr}).  The {\em symplectic
normal bundle} of $L$ in $(M,\alpha)$ is the vector bundle $$
\nu(L):=\bigcup_{x \in L} (T_xL^{d \alpha})/T_xL.
$$
It is naturally a symplectic vector bundle.
\end{definition}

\begin{remark}
If a group $G$ acts on a manifold $M$ preserving a contact form
$\alpha$ and $L$ is an invariant, isotropic submanifold,
 then the symplectic normal
bundle $\nu(L)$ has a natural $G$ action which preserves the
symplectic structure on $\nu (L)$.
\end{remark}

In \cite{wein} Weinstein proved that the symplectic normal bundle
completely determines the contact germ of the corresponding isotropic
embedding.  Theorem~\ref{contact_normal_form} is an equivariant
version of Weinstein's result in the special case where the isotropic
submanifold is a group orbit.  We now proceed with the details of the
proof of Theorem~\ref{contact_normal_form}.

The proof is a slight modification of the standard argument in the
symplectic case.  We use the notation in the statement of the theorem.
Let $Y$ denote the Reeb vector field.  Choose a $G$-invariant complex
structure $J$ on the symplectic vector bundle $\ker \alpha \to M$
compatible with the symplectic structure, that is, choose it so that
$d\alpha (\cdot, J \cdot)|_{\ker \alpha }$ is a positive definite
inner product.  The inner product is $G$-invariant by construction.
By declaring $Y$ to be a unit vector field orthogonal to $\ker \alpha$
we can extend the inner product on $\ker \alpha$ to a $G$-invariant
metric $g$ on $M$.

Because the point $x$ lies in the zero level set of the moment map, the
orbit $G\cdot x$ is an isotropic submanifold of $(M, \alpha)$.  Since
$T_x (G\cdot x)$ is isotropic in $(\ker
\alpha _x, d\alpha _x)$, $T_x (G\cdot x)$ and $J(T_x (G\cdot x))$ 
are $g$-orthogonal by construction of $g$.  Moreover 
$T_x (G\cdot x)\oplus J(T_x (G\cdot x))$ is a symplectic subspace of $(\ker
\alpha _x, d\alpha _x)$, and the $g$-orthogonal complement $V$ to 
$T_x (G\cdot x)\oplus J(T_x (G\cdot x))$ is also a symplectic
subspace.  The isotropy group $H$ of $x$ acts on $\ker \alpha _x$
preserving the decomposition $ \ker \alpha _x = T_x (G\cdot x)\oplus
J(T_x (G\cdot x)) \oplus V$ and the symplectic form $d\alpha _x
|_{\ker \alpha _x}$.  Note that $V$ is isomorphic to $(T_x (G\cdot
x))^{d\alpha}/T_x (G\cdot x)$.  Consequently $G\times _H V \to G\cdot
x = G/H$ is the symplectic normal bundle of the embedding $G\cdot x
\hookrightarrow (M, \alpha)$ while 
$G\times _H (J(T_x (G\cdot x) \times V \times \R Y_x)$ is the
topological normal bundle of the embedding $G\cdot x \hookrightarrow
M$.  The map $v\mapsto d\alpha _x (v, \cdot)|_{T_x (G\cdot x)}$
identifies $J(T_x (G\cdot x))$ with $T^*_x (G\cdot x) \simeq
(\fg/\fh)^* \simeq \fh^\circ$, where as before $\fh^\circ$ denotes the
annihilator of $\fh$ in $\fg^*$.  Therefore the topological normal
bundle of the embedding of the orbit $G\cdot x$ in $M$ is $G\times _H
(\fh^\circ \times V\times \R)$.

Denote the restriction of $d\alpha _x$ to $V$ by $\omega _V$.  Let $(W
,\beta)$ denote the contactization of $(V, \omega_V)$: $W = V\times
\R$ and $\beta = \frac{1}{2}\iota (R) \omega _V + dt$ where $R$ is the
radial vector field on $V$ and $t$ is the variable along $\R$. Thus
the normal bundle of the embedding $G\cdot x\hookrightarrow M$ can be
written as $G\times _H (\fh^\circ \times W)$.

Remark that our various identifications are chosen in such a way that 
$$
\ker \alpha _x = T_x (G\cdot x) \oplus \fh^\circ \oplus V
$$ 
and 
$$ 
d\alpha _x |_{\ker \alpha _x} = \omega + \omega _V, 
$$ 
where $\omega$ is the standard symplectic form on $T_x (G\cdot x) \oplus
\fh^\circ = (\fg/\fh) \times (\fg/\fh)^*$.

The exponential map defined by the metric $g$ provides us with a
$G$-equivariant diffeomorphism $\sigma$ from a neighborhood $U_0$ of
the zero section in $\Y = G\times _H (\fh^\circ \times W)$  to a
neighborhood $U$ of $G\cdot x$ in $M$ such that $\sigma ([1, 0, 0]) =
x$ and $d\sigma _{[1, 0, 0]} : T_x (G\cdot x) \oplus \fh^\circ \oplus
W \to T_x M$ is the identification above.

Therefore, by Darboux theorem (Theorem~\ref{ercd} above), in order to
finish proving Theorem~\ref{contact_normal_form}, it is enough to
construct on $\Y$ a $G$-invariant contact form $\epsilon$ with the
following properties:
\begin{enumerate}
\item $\ker \epsilon_{[1, 0, 0]} = T_{[1, 0, 0]} (G/H) \oplus \fh^\circ 
\oplus V$ where we identified $T_{[1, 0, 0]} \Y $ with 
$T_{[1, 0, 0]} (G/H) \oplus \fh^\circ \oplus W = \fg/\fh \oplus
\fh^\circ \oplus W$ (we think of $G/H$ as embedded in $\Y$ as the zero
section of $\Y \to G/H$).

\item $d \epsilon _{[1, 0, 0]} |_{\ker \epsilon_{[1, 0, 0]}} = 
\omega + \omega _V$

\item The moment map $F$ for the action of $G$ on $(\Y, \epsilon)$ is 
$F([g, \eta, w]) = Ad^\dagger (g)\left(\eta + i (\Phi_W (w))\right)$
where as before $ \Phi_W$ is the moment map for the action of $H$ on
the contact vector space $W$.
\end{enumerate}

The construction of $\epsilon $ is standard.  The cotangent bundle
$T^*G$ is an exact symplectic manifold with a $G\times H$-invariant
symplectic potential, where $G$ acts by the lift of the left
multiplication and $H$ by the lift of the multiplication on the right
by the inverse.  The group $H$ also acts contactly on the contact
vector space $W$.  Hence the diagonal action of $H$ on $T^*G \times W$
is contact.  It is also free.  It commutes with the trivial extension
of the action of $G$ on $T^*G$ to $T^*G \times W$.  Therefore the
contact quotient $(T^*G \times W)/\!/H$ is a contact manifold.
Moreover the contact form on $(T^*G \times W)/\!/H$ is $G$-invariant.
It is routine to verify that a choice of an $H$-invariant splitting
$\fg^*=\fh^\circ \oplus \fh^*$ allows one to identify the quotient
$(T^*G \times W)/\!/H$ with $G\times _H (\fh^\circ \times W) = \Y$ and
the induced contact form $\epsilon $ on $\Y$ has the desired
properties.  This finishes the proof of
Theorem~\ref{contact_normal_form}.


\begin{thebibliography}{WWWW}

\bibitem[Al]{?} C. Albert, Le th\'eor\`eme de r\'eduction de
Marsden-Weinstein en g\'eom\'etrie cosymplectique et de contact, {\em
J. Geom. Phys.} {\bf 6} (1989), no. 4, 627--649.


\bibitem[ACG]{ACG} J. Arms, R. Cushman and M. Gotay,  
A universal reduction procedure for Hamiltonian group actions in {\bf
The geometry of Hamiltonian systems} (Berkeley, CA, 1989), 33--51,
Math.\ Sci.\ Res.\  Inst.\ Publ., 22, Springer, New York, 1991.
 
\bibitem[AMM]{AMM} J. Arms, J. Marsden and V. Moncrief, 
Symmetry and bifurcations of momentum mappings {\em  
Comm.\ Math.\ Phys.} {\bf 78} (1980/81), no. 4, 455--478.

\bibitem[BL]{BL} L. Bates and E. Lerman, Proper group actions and 
symplectic stratified spaces, {\em Pac.\ J.\ of Math.\ (2)} {\bf 181}
(1997), 201--229.

\bibitem[Ge]{Geiges} H. Geiges, Constructions of contact manifolds,
{\em Math.\ Proc.\ Cambridge Philos.\ Soc.} {\bf 121} (1997), no. 3,
455--464.

\bibitem[GM]{GM} M. Goresky and R. MacPherson, {\em Stratified 
Morse Theory}, Berlin; New York: Springer-Verlag, 1988.

\bibitem[GS]{GS} V. Guillemin and S. Sternberg, Homogeneous quantization 
and multiplicities of group representations, {\em 
J. Funct.\ Anal.} {\bf 47} (1982), no. 3, 344--380.

\bibitem[K1]{Kbook} F. Kirwan, {\em Cohomology of Quotients in Symplectic
and Algebraic Geometry,} Mathematical Notes {\bf 31}, Princeton
University Press, Princeton, 1984.

\bibitem[K2]{Kconv} F. Kirwan, Convexity Properties of the Moment 
Mapping~III, {\em Invent. Math.} {\bf 77} (1984), 547--552.

\bibitem[L]{IJL} E. Lerman, Contact Cuts, {\em Israel J.\ Math}, to appear.
 \hfill 
See {\tt http://xxx.lanl.gov/abs/math.SG/0002041}.\hfill 

\bibitem[O]{O} J.-P. Ortega, thesis, University of California at Santa Cruz, 
1998.

\bibitem[OR]{OR} J.-P. Ortega and T. Ratiu, in preparation.


\bibitem[P]{Palais} R. Palais, On the existence of slices for actions of 
non-compact Lie groups, {\em Ann.\ Math.} {\bf 73} (1961), no. 2,
295--322.

\bibitem[SL]{SL} R. Sjamaar and E. Lerman, Stratified symplectic 
spaces and reduction, {\em Ann.\ of Math.\ } {\bf 134} (1991), 375--442.

\bibitem[We]{wein}A. Weinstein, Contact surgery and symplectic 
handlebodies,{\em Hokkaido Math. J.}, {\bf 20} (1991), no. 2,
241--251.



\end{thebibliography}
\end{document}